\documentclass[12pt]{article}
\usepackage{amssymb}

\usepackage{amsfonts,amssymb,eucal,amsmath}
\begin{document}
\begin{center}
{\large\bf Large Even Number Represent The Sum Of Odd Primes}
\end{center}
\begin{center}
 {\ Jin Li}
 \end{center}
 \begin{center}
{College of Mathematics,Sichuan University,Chengdu, PRC}
\end{center}

\begin{center}
{\bf Abstract}
\end{center}

In this paper, I proved that
$$N=p_1+p_2+2p_3, p_1\sim N/2, p_2\sim N/2, p_3=o(N),$$
where $N$ is a large even number, and $p_i\ (i=1,2,3)$ are odd primes.

\section{\large Introduction}
\setcounter{equation}{0}

In 1930, Shnirel' man proved that every integer greater than one is the sum of a bounded number of primes. This is the first significant result on the Goldbach conjecture. In 1937, I. M. Vinogradov proved that every large odd integer is the sum of three primes. Now I prove an analogue result on even number by Vinogradov's method.

\section{\large Notation}
\setcounter{equation}{0}
\par
$p_i, p$-- prime number.

$\varepsilon$-- any sufficient small positive constant.

$N$-- a sufficiently large even number.

$a,q,r,n$-- positive integers.

$\alpha ,\beta ,t$--real variables.

$c,c_i$-- some positive constant.

$\lambda$-- a suitably choosed positive number.

$$A=N\cdot e^{-\varepsilon\sqrt{\log N}}, Q=\log^{\lambda}N,\ \ \ \ \tau =A^2N^{-1}Q^{-1}.$$

$\varphi (q)$-- the Euler function

$\mu(q)$-- the M$\ddot{\rm o}$bius function

$e(x)={\rm exp}\{2\pi ix\}.$

$C_q(m)$--the Ramanujan sum
$$\sum\limits^q_{a=1,(a,q)=1}e\left(\frac{ma}{q}\right).$$

$f=O(g)$ or $f\ll g$ or $g\gg f$, if there exists a constant $c>0$ s.t. $|f(x)|\leq c g(x)$ for all $x$ in the domain of $f$. We write $f=o(g)$, if $\lim\limits_{x\rightarrow\infty}\frac{f(x)}{g(x)}=0$.

\section{Main result and proof}
\subsection{Main theorem}

{\bf Theorem 1}\ \ There exists an arithmetic function $D(N)$ and positive constant $c_1, c_2$ and $c$ such that $c_1<D(N)<c_2$ for all sufficiently large even integers $N$, and
$$R(N,A)=2D(N)A^2+O\left(\frac{A^2}{\log^cN}\right),$$
where $R(N,A)=\sum\limits_{\tiny{
\begin{array}{l}
p_1+p_2+2p_3=N,\\
\frac{N}{2}-A<p_1,p_2\leq\frac{N}{2}+A\\
2<p_3\leq A\end{array}}}\log p_1\cdot\log p_2\cdot\log p_3$

\subsection{The singular series}

We begin by studying the arithmetic functions
\begin{eqnarray*}
&&D(N)=:\sum\limits^{\infty}_{q=1,q{\rm \,is \,odd}}\frac{\mu(q)C_q(-N)}{\varphi^3(q)}-\sum\limits^{\infty}_{q=1, q {\rm \,is \,even}}\frac{\mu(q)C_q(-N)}{\varphi^3(q)}, \\ &&G(N)=:\sum\limits^{\infty}_{q=1}\frac{\mu(q)C_q(-N)}{\varphi^3(q)}.
\end{eqnarray*}
The functions $D(N), G(N)$ are called the singular series.

{\bf Theorem 2}\ \ The singular series $D(N), G(N)$ converge absolutely and uniformly in $N$ and have the Euler product
\begin{eqnarray*}
&&D(N)=2\prod\limits_{p|N,p {\rm \,is \,odd}}\left(1-\frac{1}{(p-1)^2}\right)\cdot\prod\limits_{p\nmid N}\left(1+\frac{1}{(p-1)^3}\right).\\
&&G(N)=\prod\limits_{p|N}\left(1-\frac{1}{(p-1)^2}\right)\prod\limits_{p\nmid N}\left(1+\frac{1}{(p-1)^3}\right)
\end{eqnarray*}
Moreover, for any $q>0$.
\begin{eqnarray*}
D(N,Q)&=&:\sum\limits^Q_{q=1,q {\rm \,is \,odd}}\frac{\mu(q)C_q(-N)}{\varphi^3(q)}-\sum\limits^Q_{q=1, q {\rm \,is \,oven}}\frac{\mu(q)C_q(-N)}{\varphi^3(q)}\\
&=&D(N)+O\left(\frac{1}{Q^{1-\varepsilon}}\right).
\end{eqnarray*}

{\bf Proof}\ \ Clearly,
$$|C_q(-N)|\leq\varphi (q),$$
and so
$$\left|\frac{\mu(q)C_q(-N)}{\varphi^3(q)}\right|\leq\frac{1}{\varphi^2(q)}\ll\frac{\log^2\log q}{q^2}.$$
Thus the singular series converges absolutely and uniformly in $N$. Moreover, $$D(N)-D(N,Q)\ll\sum\limits_{q>Q}\frac{1}{\varphi^2(q)}\ll\sum\limits_{q>Q}\frac{1}{q^{2-\varepsilon}}\ll\frac{1}{Q^{1-\varepsilon}}.$$
For $C_q(-N)$ is a multiplicative function of $q$ and $C_p(-N)=\left\{\begin{array}{ll}
p-1,&p|N\\
-1,&p\nmid N
\end{array}\right.$, By Euler product \begin{eqnarray*}
G(N)&=&\prod\limits_p\left(1+\sum\limits^{\infty}_{j=1}\frac{\mu(p^j)C_{p^{j}}(-N)}{\varphi^3(p^j)}\right)\\
&=&\prod\limits_p\left(1-\frac{C_p(-N)}{\varphi^3(p)}\right)=\prod\limits_{p|N}\left(1-
\frac{1}{(p-1)^2}\right)\cdot\prod\limits_{p\nmid N}\left(1+\frac{1}{(p-1)^3}\right),
\end{eqnarray*}
since $N$ is even integer $p=2|N$, then $G(N)=0$, and
$$D(N)=2\sum\limits^{\infty}_{q=1,q {\rm \,is \, odd}}\frac{\mu(q)C_q(-N)}{\varphi^3(q)}=2\prod\limits_{p|N, p {\rm \,is \,odd}}\left(1-\frac{1}{(p-1)^2}\right)\cdot\prod\limits_{p\nmid N}\left(1+\frac{1}{(p-1)^3}\right).$$
So there exist positive constant $c_1,c_2$ such that $c_1<D(N)<c_2$.

\subsection{Decomposition into major and minor arcs.}

For $1\leq q\leq Q, 0\leq a\leq q-1$ and $(a,q)=1$ the Major arc $\mathcal{M}(q,a)$ is the interval consisting of all real numbers $\alpha\in\left[-\frac{1}{\tau},1-\frac{1}{\tau}\right]$ such that
$$\left|\alpha -\frac{a}{q}\right|\leq\frac{1}{\tau}.$$
The set of major arcs is $$\mathcal{M}=\cup^Q_{q=1}\cup^{q-1}_{a=0,(a,q)=1}\mathcal{M}(q,a)\subseteq\left[-\frac{1}{\tau},
1-\frac{1}{\tau}\right],$$
and the set of minor arcs is
$$m=\left[-\frac{1}{\tau},1-\frac{1}{\tau}\right]\backslash\mathcal{M}.$$

Let
$$F(\alpha ,A)=\sum\limits_{\frac{N}{2}-A<p\leq\frac{N}{2}+A}(\log p)e(\alpha\cdot p),$$
$$\tilde{F}(\alpha ,A)=\sum\limits_{2<p\leq A}(\log p)e(\alpha\cdot 2p).$$
By the circle method, we can express $R(N,A)$ as the integal of a trigonometric polynomial over the major and minor arcs.
\begin{eqnarray*}
R(N,A)&=&\int^{1-\frac{1}{\tau}}_{-\frac{1}{\tau}}F^2(\alpha ,A)\tilde{F}(\alpha ,A)e(-N\alpha )d\alpha\\
&=&\int_{\mathcal{M}}F^2(\alpha ,A)\tilde{F}(\alpha ,A)e(-N\alpha )d\alpha +\\
&&\int_{m}F^2(\alpha ,A)\tilde{F}(\alpha ,A)e(-N\alpha )d\alpha .
\end{eqnarray*}

\subsection{The integal over the major arcs}

{\bf Theorem 3}\ \ Let
$$u(\beta ,A)=\sum\limits_{\frac{N}{2}-A<n\leq\frac{N}{2}+A}e(\beta n), \tilde{u}(\beta ,A)=\sum\limits_{2<n\leq A}e(\beta\cdot 2n).$$
Then
\begin{eqnarray*}
J(N,A)&=&\int^{1-\frac{1}{\tau}}_{-\frac{1}{\tau}}u^2(\beta ,A)\tilde{u}(\beta ,A)e(-N\beta )d\beta\\
&=&2A^2+O(A).
\end{eqnarray*}

{\bf Proof}\ \ By circle method,
\begin{eqnarray*}
J(N,A)&=&\sum\limits_{\tiny{
\begin{array}{l}
N=n_1+n_2+2n_3\\
\frac{N}{2}-A<n_1,n_2\leq\frac{N}{2}+A\\
2<n_3\leq A
\end{array}}}1\\
&=&\sum\limits_{\frac{N}{2}-A<n_1\leq\frac{N}{2}}\cdot\sum\limits_{\tiny
\begin{array}{l}
n_2+2n_3=N-n_1\\
2<n_3\leq A
\end{array}}1+\sum\limits_{\frac{N}{2}<n_1\leq\frac{N}{2}+A}\cdot\sum\limits_{\tiny
\begin{array}{l}
n_2+2n_3=N-n_1\\
2<n_3\leq A\end{array}}1\\
&=&\sum\limits_{0\leq d_1<A}\sum\limits_{\tiny\begin{array}{l}
d_2+2n_3=d_1\\
-A<d_2\leq A\\
2<n_3\leq A\end{array}}1+\sum\limits_{0<d_1\leq A}\sum\limits_{\tiny\begin{array}{l}
d_2+2n_3=-d_1\\
-A<d_2\leq A\\
2<n_3\leq A\end{array}}1\\
&=&\sum\limits_{0\leq d_1<A}(A+d_1)+\sum\limits_{0<d_1\leq A}(A-d_3)+O(A)\\
&=&2A^2+O(A).
\end{eqnarray*}

{\bf Theorem 4}(Siegel-Walfisz)

If $1\leq q\leq Q=\log^{\lambda}N$ and $(q,r)=1$, then there exists a constant $c_3$ depends only on $\lambda$ such that.
$$\theta(t;q,r)=:\sum\limits_{\tiny\begin{array}{l}
p\leq t\\
p\equiv r({\rm mod}\ q)
\end{array}}\log p=\frac{t}{\varphi (q)}+r(t),\ \ \ \ t\geq 2.$$
where $r(t)\ll t\cdot e^{-c_3\sqrt{\log t}}$.

{\bf Lemma 5}\ \ If $\alpha\in\mathcal{M}(q,a)$ and $\beta =\alpha -\frac{a}{q}$, then
$$F(\alpha ,A)=\frac{\mu (q)}{\varphi (q)}\sum\limits_{\frac{N}{2}-A<n\leq\frac{N}{2}+A}e(\beta n)+O(A\cdot e^{-c_4\sqrt{\log N}})$$

$$\tilde{F}(\alpha ,A)=\frac{C_q(2a)}{\varphi (q)}\sum\limits_{2<n\leq A}e(\beta\cdot 2n)+O(A\cdot e^{-c_4\sqrt{\log N}}).$$

{\bf Proof}\ \ Let $p\equiv r({\rm mod}\ q)$. Then $p|q$ is and only if $(r,q)>1$, and so
$$\sum\limits^q_{\tiny\begin{array}{l}
r=1\\
(r,q)=1\end{array}}\sum\limits_{\tiny\begin{array}{l}
\frac{N}{2}-A<p\leq\frac{N}{2}+A\\
p\equiv r({\rm mod}\ q)
\end{array}}e(\alpha p)\cdot\log p=\sum\limits_{\tiny\begin{array}{l}
\frac{N}{2}-A<p\leq\frac{N}{2}+A\\
p|q\end{array}}e(\alpha p)\cdot\log p\ll\sum\limits_{p|q}\log p\ll\log q.$$
So
\begin{eqnarray*}
F(\alpha ,A)&=&\sum\limits^q_{\tiny\begin{array}{l}
r=1\\
(r,q)=1\end{array}}\cdot\sum\limits_{\tiny\begin{array}{l}
\frac{N}{2}-A<p\leq\frac{N}{2}+A\\
p\equiv r({\rm mod}\ q)\end{array}}e\left(\left(\frac{a}{q}+\beta\right)\cdot p\right)\cdot\log p+O(\log q)\\
&=&\sum\limits^q_{\tiny\begin{array}{l}
r=1\\
(r,q)=1\end{array}}e(\frac{a}{q}r)\sum\limits_{\tiny\begin{array}{l}
\frac{N}{2}-A<p\leq\frac{N}{2}+A\\
p\equiv r({\rm mod}\ q)\end{array}}e(\beta p)\log p+O(\log q)\\
&=&\sum\limits^q_{\tiny\begin{array}{l}
r=1\\
(r,q)=1\end{array}}e(\frac{a}{q}r)\int^{\frac{N}{2}+A}_{\frac{N}{2}-A}e(\beta t)d\theta (t;q,r)+O(\log q)
\end{eqnarray*}
By theorem 4,
$$F(\alpha ,A)=\frac{\mu (q)}{\varphi (q)}\int^{\frac{N}{2}+A}_{\frac{N}{2}-A}e(\beta t)dt+O(Ae^{-c_3\sqrt{\log N}}).$$

\begin{eqnarray*}
\int^{\frac{N}{2}+A}_{\frac{N}{2}-A}e(\beta t)dt-\sum\limits_{\frac{N}{2}-A<n\leq\frac{N}{2}+A}e(\beta n)&=&\sum\limits_{\frac{N}{2}-A<n\leq\frac{N}{2}+A}\int^{n+1}_n(e(\beta t)-e(\beta n))dt+O(1)\\
&=&\sum\limits_{\frac{N}{2}-A<n\leq\frac{N}{2}+A}\int^{n+1}_n\left(\int^t_nde(\beta u)\right)dt+O(1)\\
&=&\sum\limits_{\frac{N}{2}-A<n\leq\frac{N}{2}+A}\int^{n+1}_n\left(\int^t_n2\pi i\beta e(\beta u)du\right)dt+O(1)\\
&\ll&\sum\limits_{\frac{N}{2}-A<n\leq\frac{N}{2}+A}|\beta |+O(1),
\end{eqnarray*}
since
$$|\beta |\leq\frac{1}{\tau}\ll A\cdot\frac{1}{\tau}=\frac{NQ}{A}\ll e^{c_4\sqrt{\log N}};$$
Similarly, we have
$$\tilde{F}(\alpha ,A)=\frac{c_q(2a)}{\varphi (q)}\sum\limits_{2<n\leq A}e(\beta\cdot 2n)+O(A\cdot e^{-c_4\sqrt{\log N}}).$$

{\bf Theorem 6}\ \ There is a positive constant $0<\varepsilon <1$, the integal over the major arcs is
$$\int_{\mathcal{M}}F^2(\alpha ,A)\tilde{F}(\alpha ,A)e(-N\alpha )d\alpha =2D(N)A^2+O\left(\frac{A^2}{\log^{(1-\varepsilon )\lambda}N}\right).$$

{\bf Proof}
\begin{eqnarray*}
&&\int_{\mathcal{M}}[F^2(\alpha ,A)\tilde{F}(\alpha ,A)-\frac{\mu^2(q)}{\varphi^3(q)}C_q(2a)u^2(\beta ,A)\tilde{u}(\beta ,A)]e(-N\alpha )d\alpha\\
&=&\sum\limits^Q_{q=1}\sum\limits^{q-1}_{\tiny\begin{array}{l}
a=0\\
(a,q)=1\end{array}}\int_{\mathcal{M}(q,a)}[F^2(\alpha ,A)\tilde{F}(\alpha ,A)-\frac{\mu^2(q)}{\varphi^3(q)}C_q(2a)u^2(\beta ,A)\tilde{u}(\beta ,A)]e(-N\alpha )d\alpha\\
&\ll&\sum\limits^Q_{q=1}\sum\limits^{q-1}_{\tiny\begin{array}{l}
a=0\\
(a,q)=1\end{array}}\int_{\mathcal{M}(q,a)}A^3\cdot e^{-c_5\sqrt{\log N}}d\alpha\\
&\leq&Q^2A^3e^{-c_5\sqrt{\log N}}\cdot\frac{1}{\tau}\\
&=&NAQ^3e^{-c_5\sqrt{\log N}}\\
&\leq&A^2e^{-c_6\sqrt{\log N}}.
\end{eqnarray*}

\begin{eqnarray*}
&&\int_{\mathcal{M}}\frac{\mu^2(q)}{\varphi^3(q)}C_q(2a)u^2(\beta ,A)\tilde{u}(\beta ,A)e(-N\alpha )d\alpha\\
&=&\sum\limits^Q_{q=1}\sum\limits^{q-1}_{\tiny\begin{array}{l}
a=0\\
(a,q)=1\end{array}}\frac{\mu^2(q)C_q(2a)}{\varphi^2(q)}\int^{\frac{a}{q}+\frac{1}{\tau}}_{\frac{a}{q}-\frac{1}{\tau}}u^2(\alpha -\frac{a}{q},A)\tilde{u}(\alpha -\frac{a}{q},A)e(-N\alpha )d\alpha\\
&=&\sum\limits^Q_{\tiny\begin{array}{l}
q=1\\
q\ {\rm \,is\, odd}\end{array}}\frac{\mu (q)}{\varphi^3(q)}\sum\limits^{q-1}_{\tiny\begin{array}{l}
a=0\\
(a,q)=1\end{array}}e(-\frac{a}{q}N)\int^{\frac{1}{\tau}}_{-\frac{1}{\tau}}u^2(\beta ,A)\\
&&\tilde{u}(\beta ,A)e(-N\beta )d\beta -\sum\limits^Q_{\tiny\begin{array}{l}
q=1\\
q\ {\rm \,is \,even}\end{array}}\frac{\mu (q)}{\varphi^3(q)}\sum\limits^{q-1}_{\tiny\begin{array}{l}
a=0\\
(a,q)=1\end{array}}e(-\frac{a}{q}N)\\
&&\int^{\frac{1}{\tau}}_{-\frac{1}{\tau}}u^2(\beta ,A)\tilde{u}(\beta ,A)e(-N\beta )d\beta\\
&=&D(N,Q)\int^{\frac{1}{\tau}}_{-\frac{1}{\tau}}u^2(\beta ,A)\tilde{u}(\beta ,A)e(-N\beta )d\beta ,
\end{eqnarray*}
because of
$$C_q(a)=\mu\left(\frac{q}{(a,q)}\right)\varphi (q)\varphi^{-1}\left(\frac{q}{(a,q)}\right).$$

If $|\beta |\leq\frac{1}{2}$,then
$$u(\beta ,A)\ll\frac{1}{|\beta |}, \tilde{u}(\beta ,A)\ll\frac{1}{|\beta |}$$
and
$$\int^{\frac{1}{2}}_{\frac{1}{\tau}}u^2(\beta ,A)\tilde{u}(\beta ,A)e(-N\beta )d\beta\ll\int^{\frac{1}{\tau}}_{\frac{1}{\tau}}\frac{d\beta}{\beta^3}\ll\tau^2=A^2Q^{-2}e^{-2\varepsilon\sqrt{\log N}}.$$
Similarly,
$$\int^{-\frac{1}{\tau}}_{-\frac{1}{2}}u^2(\beta ,A)\tilde{u}(\beta ,A)e(-N\beta )d\alpha\ll A^2Q^{-2}e^{-2\varepsilon\sqrt{\log N}}.$$
By theorem 3,
\begin{eqnarray*}
&&\int^{\frac{1}{\tau}}_{-\frac{1}{\tau}}u^2(\beta ,A)\tilde{u}(\beta ,A)e(-N\beta )d\beta \\
&=&\int^{\frac{1}{2}}_{-\frac{1}{2}}u^2(\beta ,A)\tilde{u}(\beta ,A)e(-N\beta )d\beta +O(A^2e^{-2\varepsilon\sqrt{\log N}})\\
&=&2A^2+O(A^2e^{-2\varepsilon\sqrt{\log N}}).
\end{eqnarray*}
By theorem 2,
$$D(N,Q)=D(N)+O\left(\frac{1}{Q^{1-\varepsilon}}\right).$$
Therefore,
\begin{eqnarray*}
&&\int_{\mathcal{M}}F^2(\alpha ,A)\tilde{F}(\alpha ,A)e(-N\alpha )d\alpha\\
&=&D(N,Q)\int^{\frac{1}{\tau}}_{-\frac{1}{\tau}}u^2(\beta ,A)\tilde{u}(\beta ,A)e(-N\beta )d\beta +O(A^2e^{-c_6\sqrt{\log N}})\\
&=&2D(N)A^2+O(\frac{A^2}{\log^{(1-\varepsilon )\lambda}N}).
\end{eqnarray*}
This completes the proof.

\subsection{Proof of theorem}

{\bf Theorem 7}(I. M. Vinogradov)\ \ If $\left|\alpha -\frac{a}{q}\right|\leq\frac{1}{q^2}$, where $a$ and $q$ are integers such that $1\leq q\leq A$ and $(a,q)=1$, then
$$\tilde{F}(\alpha ,A)=\sum\limits_{2<p\leq A}e(\alpha\cdot 2p)\log p\ll A\log^4N\left(\sqrt{\frac{q}{A}}+\sqrt{\frac{1}{q}}+\frac{1}{H}\right),$$
where
$$H=e^{\frac{1}{2}\sqrt{\log N}}.$$

{\bf Lemma 8}(Dirichlet)\ \ If $\alpha\in m$, then there must exist integer $q,a$ such that
$$\frac{a}{q}\in\left[-\frac{1}{\tau},1-\frac{1}{\tau}\right], (a,q)=1, Q<q\leq\tau$$
and
$$\left|\alpha-\frac{a}{q}\right|<\frac{1}{q\tau}.$$

{\bf Theorem 9}\ \ For any $\lambda >0$, we have
$$\int_{m}F^2(\alpha ,A)\tilde{F}(\alpha ,A)e(-\alpha N)d\alpha\ll\frac{A^2}{\log^{\frac{\lambda}{2}-5}}.$$

{\bf Proof}\ \ If $\alpha\in m$, then $Q<q\leq\tau$. By theorem 7,
\begin{eqnarray*}
\tilde{F}(\alpha ,A)&\ll&A\log^4N\left(\sqrt{\frac{q}{A}}+\sqrt{\frac{1}{q}}+\frac{1}{H}\right)\ll A\log^4N\\
\left(\sqrt{\frac{\tau}{A}}+\sqrt{\frac{1}{Q}}+\frac{1}{H}\right)&=&A\log^4N\left(e^{-\frac{\varepsilon}{2}\sqrt{\log N}}\log^{-\frac{\lambda}{2}}N+\log^{-\frac{\lambda}{2}}N+e^{-\frac{1}{2}\sqrt{\log N}}\right)\\
&\ll&\frac{A}{\log^{\frac{\lambda}{2}-4}N}.
\end{eqnarray*}
Since
$$\theta(N,A)=:\sum\limits_{\frac{N}{2}-A<p\leq\frac{N}{2}+A}\log p\ll A,$$
we have
$$\int^{1-\frac{1}{\tau}}_{-\frac{1}{\tau}}|F(\alpha ,A)|^2d\alpha =\sum\limits_{\frac{N}{2}-A<p\leq\frac{N}{2}+A}\log^2p\ll\log N\cdot\sum\limits_{\frac{N}{2}-A<p\leq\frac{N}{2}+A}\log p\ll A\log N,$$
and so
\begin{eqnarray*}
\int_m|F(\alpha ,A)|^2|\tilde{F}(\alpha ,A)|d\alpha&\ll&\sup\{|\tilde{F}(\alpha ,A)|:\alpha\in m\}\int_m|F(\alpha ,A)|^2d\alpha\\
&\ll&\frac{A}{\log^{\frac{\lambda}{2}-4}N}\int^1_0|F(\alpha ,A)|^2d\alpha\\
&\ll&\frac{A^2}{\log^{\frac{\lambda}{2}-5}N}
\end{eqnarray*}
This completes the proof.

{\bf Proof of theorem 1}\ \
\begin{eqnarray*}
R(N,A)&=&\int_{\mathcal{M}}F^2(\alpha ,A)\tilde{F}(\alpha ,A)e(-N\alpha )d\alpha +\int_mF^2(\alpha ,A)\tilde{F}(\alpha ,A)e(-N\alpha)d\alpha\\
&=&2D(N)A^2+O\left(\frac{A^2}{\log^{(1-\varepsilon)\lambda}N}\right)+O\left(\frac{A^2}{\log^{\frac{\lambda}{2}-5}}\right),
\end{eqnarray*}
let
$$c=min\{(1-\varepsilon )\lambda ,\frac{\lambda}{2}-5\}, \lambda >10.$$

{\bf Acknowledgment.}\ \ I am grateful to my parents. I also wish to thank professor shiqing zhang, from whom I learned a lot.

\end{document}